\documentclass[graybox]{svmult}

\usepackage{type1cm}          
\usepackage{makeidx}          
\usepackage{graphicx}         
\usepackage{multicol}         
\usepackage[bottom]{footmisc} 

\usepackage{newtxtext}     
\usepackage[varvw]{newtxmath} 

\usepackage{xparse}
\usepackage{ifthen}
\usepackage{hyperref}
\usepackage{breakurl}

\usepackage{amsmath}
\usepackage{bm}

\usepackage{booktabs}
\usepackage{multirow}

\usepackage{caption}
\usepackage{subcaption}

\usepackage[ruled]{algorithm2e}
\usepackage{listings}

\makeindex 

\newcommand{\bracketed}[1]{\left[#1\right]}

\NewDocumentCommand{\iprod}{o m m}{
    \left(#2,\: #3\right)\IfValueT{#1}{_{#1}}
}

\NewDocumentCommand{\iprodS}{o m m}{
    \left\langle#2,\: #3\right\rangle\IfValueT{#1}{_{#1}}
}

\NewDocumentCommand{\norm}{o m}{
    \left\vert\left\vert#2\right\vert\right\vert\IfValueT{#1}{_{#1}}
}

\NewDocumentCommand{\transpT}{o m}{
    #2\IfNoValueTF{#1}{^\mathrm{T}}{^{\overset{#1}{\mathrm{T}}}}
}

\NewDocumentCommand{\tr}{o m}{
    \mathrm{tr}\IfNoValueTF{#1}{\: #2}{\bracketed{#2}}
}

\NewDocumentCommand{\sym}{o m}{
    \mathrm{sym}\IfNoValueTF{#1}{\: #2}{\bracketed{#2}}
}

\NewDocumentCommand{\skw}{o m}{
    \mathrm{skw}\IfNoValueTF{#1}{\: #2}{\bracketed{#2}}
}

\NewDocumentCommand{\as}{o m}{
    \mathrm{as}\IfNoValueTF{#1}{\: #2}{\bracketed{#2}}
}

\NewDocumentCommand{\determ}{o m}{
    \mathrm{det}\IfNoValueTF{#1}{\: #2}{\bracketed{#2}}
}

\NewDocumentCommand{\inv}{o m}{
    \IfNoValueTF{#1}{#2 ^{-1}}{\bracketed{#2}^{-1}}
}

\NewDocumentCommand{\invt}{o m}{
    \IfNoValueTF{#1}{#2 ^{-\mathrm{T}}}{\bracketed{#2}^{-\mathrm{T}}}
}

\DeclareMathOperator{\grad}{\nabla}

\renewcommand{\div}{\nabla\cdot}

\NewDocumentCommand{\divc}{o m m}{
    \nabla_{\vec{#3}} \cdot \IfNoValueTF{#1}{#2}{\bracketed{#2}}
}

\NewDocumentCommand{\gradc}{o m m}{
    \nabla_{\vec{#3}} \: \IfNoValueTF{#1}{#2}{\bracketed{#2}}
}

\NewDocumentCommand{\mtDiff}{o m}{
    \IfNoValueTF{#1}{#2'}{\left(#2\right)'}
}

\NewDocumentCommand{\VecSpace}{o o o m}{
    \IfNoValueTF{#3}
    {\mathrm{#4} \IfValueT{#2}{_{#2}} \ifthenelse{\equal{#1}{}}{}{\left(#1\right)}}
    {\left( \mathrm{#4} \IfValueT{#2}{_{#2}} \ifthenelse{\equal{#1}{}}{}{\left(#1\right)} \right) ^{#3}}
}

\NewDocumentCommand{\LiiSpace}{o o o}{
    \IfNoValueTF{#3}
    {\mathrm{L}^2 \IfValueT{#2}{_{#2}} \ifthenelse{\equal{#1}{}}{}{\left(#1\right)}}
    {\left( \mathrm{L}^2 \IfValueT{#2}{_{#2}} \ifthenelse{\equal{#1}{}}{}{\left(#1\right)} \right) ^{#3}}
}

\NewDocumentCommand{\HiSpace}{o o o}{
    \IfNoValueTF{#3}
    {\mathrm{H}^1 \IfValueT{#2}{_{#2}} \IfValueT{#1}{\left(#1\right)}}
    {\left( \mathrm{H}^1 \IfValueT{#2}{_{#2}}\IfValueT{#1}{\left(#1\right)} \right) ^{#3}}
}

\NewDocumentCommand{\HdivSpace}{o o o}{
    \IfNoValueTF{#3}
    {\mathrm{H} \IfValueT{#2}{_{#2}} \left(\mathrm{div}\ifthenelse{\equal{#1}{}}{}{,\, #1}\right)}
    {\left( \mathrm{H} \IfValueT{#2}{_{#2}} \left(\mathrm{div}\ifthenelse{\equal{#1}{}}{}{,\, #1}\right) \right) ^{#3}}
}

\NewDocumentCommand{\feP}{o o}{
    \IfNoValueTF{#2}
    {\mathrm{P} \IfValueT{#1}{_{#1}}}
    {\left( \mathrm{P} \IfValueT{#1}{_{#1}} \right)^{#2}}
}

\NewDocumentCommand{\feDP}{o o}{
    \IfNoValueTF{#2}
    {\mathrm{DP} \IfValueT{#1}{_{#1}}}
    {\left( \mathrm{DP} \IfValueT{#1}{_{#1}} \right)^{#2}}
}

\NewDocumentCommand{\feRT}{o o}{
    \IfNoValueTF{#2}
    {\mathrm{RT} \IfValueT{#1}{_{#1}}}
    {\left( \mathrm{RT} \IfValueT{#1}{_{#1}} \right)^{#2}}
}

\NewDocumentCommand{\feBDM}{o o}{
    \IfNoValueTF{#2}
    {\mathrm{BDM} \IfValueT{#1}{_{#1}}}
    {\left( \mathrm{BDM} \IfValueT{#1}{_{#1}} \right)^{#2}}
}

\renewcommand{\vec}[1]{\bm{\mathrm{#1}}}
\newcommand{\vecGreek}[1]{\bm{#1}}
\newcommand{\mat}[1]{\bm{\mathrm{#1}}}
\newcommand{\matGreek}[1]{\bm{#1}}

\newcommand{\spEqlb}{\mathrm{R}}

\newcommand{\h}{h}

\newcommand{\strainGreenLagrangeLin}{\matGreek{\varepsilon}}
\newcommand{\Lamei}{\lambda}

\newcommand{\stressPKLin}{\matGreek{\sigma}}

\newcommand{\LinPKEqlbH}{\matGreek{\sigma}^{\mathrm{R}}_\h}
\newcommand{\Displacement}{\vec{u}}

\newcommand{\domain}{\Omega}
\newcommand{\domainH}{\mathcal{T}_{\h}}

\newcommand{\patch}{\omega_{z}}
\newcommand{\hatfunc}{\varphi_{z}}

\newcommand{\GenSol}{\mathrm{u}}
\newcommand{\GenSolH}{\mathrm{u}_\h}

\newcommand{\flux}{\vecGreek{\varsigma}}
\newcommand{\fluxH}{\vecGreek{\varsigma}\left(\GenSolH\right)}
\newcommand{\fluxHi}{\vecGreek{\varsigma}\left(\GenSolH\right)\big\vert_i}

\newcommand{\fluxR}{\vecGreek{\varsigma}^\spEqlb}
\newcommand{\fluxRH}{\vecGreek{\varsigma}^\spEqlb_\h}

\newcommand{\fluxRHzSEi}{\Delta\widetilde{\bm{\varsigma}}^R_{z,h}}
\newcommand{\fluxRHzSEii}{\Delta\bm{\varsigma}^R_{z,h}}
\newcommand{\fluxRHi}{\vecGreek{\varsigma}^\spEqlb_\h\big\vert_i}
\newcommand{\fluxRHziSEi}{\Delta\widetilde{\bm{\varsigma}}^R_{z,h}\big\vert_i}
\newcommand{\fluxRHziSEii}{\Delta\bm{\varsigma}^R_{z,h}\big\vert_i}

\newcommand{\LinPKEqlbHz}{\matGreek{\sigma}^{\mathrm{R}}_{z,\h}}
\newcommand{\LinPKEqlbHzSEi}{\Delta\widetilde{\matGreek{\sigma}}^{\mathrm{R}}_{z,\h}}
\newcommand{\LinPKEqlbHzSEii}{\Delta\matGreek{\sigma}^{\mathrm{R}}_{z,\h}}
\newcommand{\LinPKEqlbHzWS}{\Delta_\mathrm{sym}\matGreek{\sigma}^{\mathrm{R}}_{z,\h}}

\SetKwFor{For}{for (}{)}{}

\definecolor{codegreen}{rgb}{0,0.6,0}
\definecolor{codegray}{rgb}{0.5,0.5,0.5}
\definecolor{codepurple}{rgb}{0.58,0,0.82}
\definecolor{backcolour}{rgb}{0.95,0.95,0.92}

\lstdefinestyle{codestyle}{
    backgroundcolor=\color{backcolour},   
    commentstyle=\color{codegreen},
    keywordstyle=\color{magenta},
    numberstyle=\tiny\color{codegray},
    stringstyle=\color{codepurple},
    basicstyle=\ttfamily\footnotesize,
    breakatwhitespace=false,         
    breaklines=true,                 
    captionpos=b,                    
    keepspaces=true,                 
    numbers=left,                    
    numbersep=5pt,                  
    showspaces=false,                
    showstringspaces=false,
    showtabs=false,                  
    tabsize=2
}

\begin{document}

\title*{Adaptive finite element methods based on flux and stress equilibration using FEniCSx}
\author{Maximilian Brodbeck\orcidID{0000-0002-1792-7318},\\ 
        Fleurianne Bertrand \orcidID{0000-0001-5111-5105} and\\
        Tim Ricken \orcidID{0000-0001-8515-5009}}
\institute{Maximilian Brodbeck, Tim Ricken \at Institute of Structural Mechanics and Dynamics, Pfaffenwaldring 27, 70569 Stuttgart, Germany \email{{brodbeck, ricken}@isd.uni-stuttgart.de} 
\and Fleurianne Bertrand \at Numerik partieller Differentialgleichungen,  Reichenhainer Straße 41, 09126 Chemnitz, Germany  \email{fleurianne.bertrand@mathematik.tu-chemnitz.de}
}

\maketitle
\vspace{-1cm}
\abstract{This contribution shows how a-posteriori error estimators based on equilibrated fluxes -- H(div) functions fulfilling the underlying conservation law -- can be implemented in FEniCSx. Therefore, dolfinx\_eqlb is introduced, its algorithmic structure is described and classical benchmarks for adaptive solution procedures for the Poisson problem and linear elasticity are presented.} 

\section{Introduction}
The accurate resolution of physical quantities in numerical simulations is of significant importance across various fields of engineering and applied sciences. Finite element methods have been extensively employed, and convergence and optimal a priori error estimates have been proven. 
Various software packages have been proposed, while the latest trend -- followed e.g. by FEniCSx \cite{FEniCSx_2023} -- focuses on abstraction, generality and automation without lousing computational efficiency. 
It is well known that on general domains, numerical solutions often lack the regularity, required to directly apply a priori estimates. 
To maintain optimal convergence, adaptive procedures based on a posteriori error estimates, combined with local mesh refinement, have been developed. 
FEniCSx is currently well-suited for handling residual-based error estimators, making it a valuable tool for many adaptive finite element methods. 
However, another important class of a posteriori error estimation involves the equilibration of flux or stress, leading to guaranteed, fully localized, and easily computable upper bounds on the error, especially in the energy norm.
Rooted in the hypercircle identity of Prager and Synge \cite{Prager_Equilibartion_1947}, the Poisson problem is discussed in works such as \cite{Braess_EqlbFluxes_2008}, \cite{Cai_SemiexplzEqlb_2012}, \cite{Ern_FluxEqlb_2015} or \cite{Bertrand_Hypercircle_2020} while applications to liner elasticity can be found e.g. in \cite{Bertrand_EqlbElast_2021}.
Even tough this method is powerful, a direct implementation into FEniCSx is challenging and currently not available in literature.  

This gap is addressed by introducing dolfinx\_eqlb, a library for the efficient computation of equilibrated fluxes and stresses.
Error estimates for the Poisson problem and linear elasticity alongside with requirements on the equilibration process are reviewed in section  \ref{sec:ReviewEE}. 
The basic algorithmic structure of the software is presented in section \ref{sec:ComputingEqlb}, while this contribution is concluded by discussing two benchmarks for adaptive finite element methods in section \ref{sec:Results}.

\section{A-posteriori error estimation based on equilibrated fluxes}
\label{sec:ReviewEE}

\subsection{Equilibration with full gradient}
Equilibration in presence of the full gradient starts from the Poisson problem in first order form
\begin{equation}
    \div\flux(\GenSol) = \mathrm{f} \;\text{in}\; \domain \quad\text{with}\quad \flux(\GenSol) := -\kappa\,\grad \GenSol \quad\text{and}\quad 
    \begin{cases}
        \GenSol = 0                    & \text{on} \; \Gamma_\mathrm{D}\\
        \flux(\GenSol)\cdot\vec{n} = 0 & \text{on} \; \Gamma_\mathrm{N}\\
    \end{cases}\; , 
    \label{eq:poisson}
\end{equation}
as it is discussed e.g. in \cite{Braess_EqlbFluxes_2008} or \cite{Ern_FluxEqlb_2015}. For any $\mathrm{f} \in \LiiSpace[\domain]$, the weak solution $\GenSol \in \HiSpace[\domain][\Gamma_\mathrm{D}]$ satisfies
\begin{equation}
    \iprod{\flux(\GenSol)}{\grad\mathrm{v}} = \iprod{\mathrm{f}}{\mathrm{v}} \quad\text{for all}\quad \mathrm{v} \in \HiSpace[\domain][\Gamma_\mathrm{D}]\; .
    \label{eq:poisson_weak}
\end{equation}
Defining the space $\VecSpace[]{V}_k := \left\{\mathrm{v}_\h \in \HiSpace[\domain]:\, \mathrm{v}_\h\vert_\mathrm{T} \in \feP[k](\mathrm{T})
\right\}$, $\GenSol$ can be approximated in $\VecSpace[]{V}_{\Gamma ,k}  := \VecSpace[]{V}_{k} \cap \HiSpace[\domain][\Gamma_\mathrm{D}]$.
For further deviations, the equilibrated flux is introduced:
\begin{definition}
    An equilibrated flux is a function $\fluxR\in\HdivSpace[\domain][\Gamma_\mathrm{N}]$ 
    fulfilling 
    \begin{equation}
        \div\fluxR = \mathrm{f} \;\text{in }\; \Omega \quad\text{and}\quad \fluxR \cdot \vec{n} = 0 \; \text{on} \; \Gamma_\mathrm{N} \ .
        \label{eq:flux_eqlb_conditions} 
    \end{equation}
    \label{def:equilibrated_flux}
\end{definition}
\vspace{-0.7cm}

\noindent For any arbitrary approximation $\GenSolH \in \VecSpace[]{V}_{\Gamma ,k}$, any equilibrated flux $\fluxR$ and the differences $\delta\flux = \flux(\GenSol) - \fluxR$ respectively $\delta\GenSol = \GenSol - \GenSolH$, the so called Prager-Synge identity
\begin{equation}
    \iprod{\delta\flux}{\grad\delta\GenSol} = \iprodS[\partial\domain]{\delta\flux\cdot\vec{n}}{\delta\GenSol} - \iprod{\div\delta\flux}{\delta\GenSol} = 0\; ,
    \label{eq:prager_synge_identity}
\end{equation}
is satisfied. Based thereon, the following error estimate holds:
\begin{theorem}
    Let $\kappa$ be constant, $\GenSol \in \HiSpace[\domain][\Gamma_\mathrm{D}]$ be the solution of \eqref{eq:poisson_weak}, $\GenSolH \in \VecSpace[]{V}_{\Gamma ,k}$ be arbitrary and $\fluxRH\in\feRT[m]$ an approximation of $\fluxR$ following Definition \ref{def:equilibrated_flux}, then
    \begin{equation}
        \norm{\grad\bracketed{\GenSol - \GenSolH}}^2 \leq \sum_{\mathrm{T}\in\domainH} \bracketed{\norm[\mathrm{T}]{\fluxRH - \fluxH} + C_P \norm[\mathrm{T}]{\mathrm{f} - \div\fluxRH}}^2 \; .
        \label{eq:ee_poisson}
    \end{equation}
    \label{thm:ee_poisson}
\end{theorem}
\vspace{-0.9cm}
\begin{remark}
    If $\kappa$ is discontinous, a modified estimate following \cite{Cai_SemiexplzEqlb_2012} holds:
    \begin{equation}
        \norm{\kappa^{1/2}\grad\bracketed{\GenSol - \GenSolH}}^2 \leq \sum_{\mathrm{T}\in\domainH} \bracketed{\norm[\mathrm{T}]{\kappa^{-1/2}\bracketed{\fluxRH - \fluxH}} + C_P \norm[\mathrm{T}]{\mathrm{f} - \div\fluxRH}}^2 \; .
        \label{eq:ee_poisson-disc-kappa}
    \end{equation}
\end{remark}

\vspace{-0.8cm}
\subsection{Equilibration with weakly symmetric gradient}
This section introduces the equilibration procedure when a symmetric gradient $\strainGreenLagrangeLin(\Displacement) = \sym{\grad\Displacement}$ is present. 
Based on the linearised Piola-Kirchoff stress $\stressPKLin(\Displacement) = 2\strainGreenLagrangeLin(\Displacement) + \tilde{\Lamei}\,\div\Displacement\,\mat{I}$,  the balance of linear momentum reads
\begin{equation}
    \div\stressPKLin(\Displacement) = -\vec{f} \;\text{in}\; \domain \quad\text{with}\quad \Displacement = \vec{0} \; \text{on} \; \Gamma_\mathrm{D}  \quad\text{and}\quad \stressPKLin(\Displacement) \cdot \vec{n} = \vec{t} \; \text{on} \; \Gamma_\mathrm{N}\; .
    \label{eq:linear_elasticity}
\end{equation}

\noindent For any $\vec{f} \in \LiiSpace[\domain]$, the weak solution $\Displacement \in \HiSpace[\domain][\Gamma_\mathrm{D}][2]$ fulfills
\begin{equation}
    \iprod{\stressPKLin(\Displacement)}{\strainGreenLagrangeLin(\vec{v})} = \iprod{\vec{f}}{\vec{v}} - 
    \iprodS[\Gamma_\mathrm{N}]{\vec{t}}{\vec{v}}
    \quad\text{for all}\quad \vec{v} \in \HiSpace[\domain][\Gamma_\mathrm{D}][2]\; .
    \label{eq:linear_elasticity-weak}
\end{equation}
Considering the symmetry of the stress tensor in a weak sense (see \cite{Bertrand_EqlbElast_2021}), the following definition of the equilibrated stress tensor is introduced:
\begin{definition}
An equilibrated stress is a function $\LinPKEqlbH\in\feRT[m][d]$ 
    fulfilling
    \begin{equation}
        \div\LinPKEqlbH = -\Pi_{m-1}\,\vec{f}\;\text{on}\;\domain \quad\text{and}\quad  \LinPKEqlbH \cdot \vec{n} = \vec{t} \; \text{on} \; \Gamma_\mathrm{N}
        \label{eq:stress_eqlb_conditions}
    \end{equation}
    and the weak symmetry condition
    \begin{equation}
        \iprod{\sigma^\spEqlb_\h\big\vert_{12} - \sigma^\spEqlb_\h\big\vert_{21}}{\gamma_\h} = 0 \quad\text{for all}\quad \gamma_\h \in \VecSpace[]{V}_{1}\; .
        \label{eq:stress_eqlb_conditions_ws}
    \end{equation}
    \label{def:equilibrated_stress}
\end{definition}
\vspace{-0.7cm}
\noindent Let further error- and $\mathcal{A}$-norm be defined by
\begin{equation}
    \vert\vert\vert \bullet \vert\vert\vert^2 = \norm{\strainGreenLagrangeLin\left(\bullet\right)}^2 + \tilde{\Lamei}\,\norm{\div\left(\bullet\right)}^2\, ,
\end{equation}
respectively
\begin{equation}
    \norm[\mathcal{A}]{\left(\bullet\right)}^2 = \iprod{\left(\bullet\right)}{\mathcal{A}\left(\bullet\right)} \;\text{with}\; \mathcal{A}\left(\bullet\right) = \frac{1}{2}\bracketed{\left(\bullet\right) - \frac{\tilde{\Lamei}}{2\,(1+\tilde{\Lamei})}\tr{\left(\bullet\right)}\,\mat{I}}\, ,
\end{equation}
the following error estimate holds:
\begin{theorem}
    Let $\Displacement \in \HiSpace[\domain][\Gamma_\mathrm{D}][2]$ be the solution of \eqref{eq:linear_elasticity-weak}, $\Displacement_\h \in \left(\VecSpace[]{V}_{\Gamma ,k}\right)^d$ be arbitrary, and $\LinPKEqlbH$ following Definition \ref{def:equilibrated_stress}. Denoting the difference $\Displacement - \Displacement_\h$ by $\delta\Displacement$, then
    \begin{equation}
        \vert\vert\vert \delta\Displacement \vert\vert\vert^2 \leq \norm[\mathcal{A}]{\LinPKEqlbH - \stressPKLin(\Displacement_\h)}^2 + C_K \sum_{\mathrm{T}\in\domainH} \bracketed{\norm[\mathrm{T}]{\skw{\LinPKEqlbH}} + C_P \norm[\mathrm{T}]{\vec{f} + \div\LinPKEqlbH}}^2\; .
        \label{eq:ee_linelasticity}
    \end{equation}
    \label{thm:ee_linear_elasticity}
\end{theorem}

\newpage
\begin{proof}
    Evaluating the $\mathcal{A}$-norm of the $\LinPKEqlbH - \stressPKLin(\Displacement_\h)$ yields
    \begin{equation}
        \norm[\mathcal{A}]{\LinPKEqlbH - \stressPKLin(\Displacement_\h)} \geq \vert\vert\vert \delta\Displacement \vert\vert\vert^2 +\norm{\strainGreenLagrangeLin\left(\delta\Displacement\right)}^2 - 2\,\iprod{\delta\stressPKLin^\spEqlb}{\strainGreenLagrangeLin\left(\delta\Displacement\right)}\; ,
        \label{eq:proof_ee-le_eval_anorm}
    \end{equation}
    where $\delta\stressPKLin^\spEqlb$ denotes the difference between true- and equilibrated stress $\stressPKLin(\Displacement) - \LinPKEqlbH$. 
    Integration by parts under consideration of the symmetry of the true stress and the equilibration conditions \eqref{eq:stress_eqlb_conditions} allows a reformulation of the mixed term:
    \begin{equation}
        \iprod{\delta\stressPKLin^\spEqlb}{\strainGreenLagrangeLin\left(\delta\Displacement\right)} = \iprod{\vec{f} + \div\LinPKEqlbH}{\delta\Displacement} + \iprod{\skw{\LinPKEqlbH}}{\grad\delta\Displacement}
        \label{eq:proof_ee-le_reformulation}
    \end{equation}
    Based on the weak symmetry \cite[Lemma 1]{Bertrand_EqlbElast_2021}, \eqref{eq:proof_ee-le_reformulation} can be bounded from above
    \begin{equation}
        \iprod{\delta\stressPKLin^\spEqlb}{\strainGreenLagrangeLin\left(\delta\Displacement\right)} \leq C_K \sum_{\mathrm{T}\in\domainH} \bracketed{\norm[\mathrm{T}]{\skw{\LinPKEqlbH}} + C_P \norm[\mathrm{T}]{\vec{f} + \div\LinPKEqlbH}}^2 + \norm{\strainGreenLagrangeLin\left(\delta\Displacement\right)}\; .
        \label{eq:proof_ee-le_bound_one}
    \end{equation}
    Inserting \eqref{eq:proof_ee-le_bound_one} into \eqref{eq:proof_ee-le_eval_anorm} completes the proof.
\end{proof}

\section{Algorithms and implementation}
\label{sec:ComputingEqlb}
This section introduces dolfinx\_eqlb, a FEniCSx based library for flux and stress equilibration. 
In order to keep the presentation general, $\bm{\theta}$ denotes followingly either a flux or a stress. Adaptive finite element methods are typically based on the loop\\
\centerline{$...$ $\rightarrow$ SOLVE $\rightarrow$ ESTIMATE $\rightarrow$ MARK $\rightarrow$ REFINE $\rightarrow$ $...$}
Using equilibration based error estimates requires
\begin{enumerate}
    \item the evaluation of projections of the right-hand-side (RHS) $\Pi_{m-1}\,\mathrm{f}$ and the approximated flux $\Pi_{m-1}\,\bm{\theta}_\h$ in a discontinuous Lagrange space of order $m-1 \geq k-1$
    \item the calculate of the equilibrated flux $\bm{\theta}^\spEqlb_\h\in\feRT[m][d]$
\end{enumerate}
during the step ESTIMATE. Therefore the constrained minimisation problem
\begin{equation}
    \bm{\theta}^\spEqlb_\h := \arg\underset{\vec{v} \in \feRT[m][d] \land \mathrm{constr}}{\min} \norm{\vec{v} - \bm{\theta}_\h}
    \label{eq:eqlb_by_global_minimisation}
\end{equation}
is solved. 
The constraints and dimension $d$ follow from Definition \ref{def:equilibrated_flux} or \ref{def:equilibrated_stress}.
As a direct solution of \eqref{eq:eqlb_by_global_minimisation} is computationally too expensive, the problem is localised by introducing for each node $z$ the nodal, piece-wise linear basis-function $\hatfunc$. 
A patch $\patch$ is then just the support of $\hatfunc$, on which a local function space
\begin{equation*}
    \mathrm{V}_m(\patch) := \left\{\vec{v}\in\feRT[m][d]:\; \vec{v}\cdot\vec{n}=
    \begin{cases}
        0                                          &  \partial\patch \cap \Gamma_\mathrm{N} = \emptyset\\
        \hatfunc\,\tilde{\mathrm{t}} & \text{else}\;
    \end{cases} 
    \right\}\, ,
\end{equation*}
is defined. 
Therein, the projection of the normal trace $\bm{\theta}\cdot\vec{n}$ into the facet-wise polynomial space of order $m-1$ is denoted by $\tilde{\mathrm{t}}$.

\newpage
\noindent Combining all $\bm{\theta}^\spEqlb_{\h,z} \in \mathrm{V}_m(\patch)$
\begin{equation}
    \bm{\theta}^\spEqlb_\h = \sum_z \bm{\theta}^\spEqlb_{\h,z} \; \text{with} \; \bm{\theta}^\spEqlb_{\h,z} := \arg\underset{\vec{v} \in \mathrm{V}_m(\patch) \land \mathrm{constr}}{\min} \norm[\patch]{\vec{v} - \hatfunc\bm{\theta}_\h}
    \label{eq:eqlb_by_local_minimisation}
\end{equation}
yields the same result as \eqref{eq:eqlb_by_global_minimisation}, replacing on large by a series of small problems.

Algorithmically, the equilibration of fluxes -- within the following, this term denotes vector-valued functions with an equilibration condition of form \eqref{eq:poisson}\textsubscript{1} -- and stresses can be handled similarly, as both share constraints on BCs and divergence.
Therefore, stresses are handled in a first step as multiple fluxes, where each 'flux' is one row of the stress tensor. 
Symmetry is enforced in a second step. 
The resulting algorithmic structure at the mesh level is described in Algorithm \ref{alg:eqlb_mesh_level}. 
Starting with lists of DOLFINx functions for the equilibrated fluxes $\left\{\fluxRHi\right\}$, the projected fluxes $\left\{\Pi_{m-1}\,\fluxHi\right\}$, the projected RHS $\left\{\Pi_{m-1}\,\mathrm{f}_i\right\}$ and the facets ($fct$) on the Dirichlet boundary of the primal problem $\left\{fct\in\Gamma_\mathrm{D}\right\}$ a general patch is created. 
It serves as sub-mesh respectively sub-function-space and has to be updated for each patch around each mesh node.
Flux equilibration and the optional enforcement of the weak symmetry condition are performed on each patch. 
The therefore required Algorithms \ref{alg:eqlb_se_patch_level} and \ref{alg:eqlb_ws_patch_level} are discussed in sections \ref{sec:Eqlb_LocMinProb} and \ref{sec:Eqlb_WeakSym}.
\vspace{-0.3cm}
\begin{algorithm}
    \caption{Equilibration on the mesh level.}
    \label{alg:eqlb_mesh_level}

    \SetKwInOut{Input}{input }
    \SetKwFunction{cPatch}{Patch}\SetKwFunction{cpatch}{patch}

    \Input{$\left\{\fluxRHi\right\}$, $\left\{\Pi_{m-1}\,\fluxHi\right\}$, $\left\{\Pi_{m-1}\,\mathrm{f}_i\right\}$ and $\left\{fct\in\Gamma_\mathrm{D}\right\}$}

    \BlankLine
    \cpatch $\gets$ \cPatch{mesh, $\left\{f\in\Gamma_\mathrm{D}\right\}$, function\_spaces}\;

    \BlankLine
    \For{$n = 0;\ n < n_\mathrm{nodes};\ i++$}
    {
        patch.create\_subdofmap($n$)\;
        equilibrate\_flux\_semiexplt(...) \tcp*{see Algorithm \ref{alg:eqlb_se_patch_level}}
        \textbf{if} \textit{weaksym\_stresses} \textbf{then} impose\_weak\_symmetry(...) \tcp*{see Algorithm \ref{alg:eqlb_ws_patch_level}}
    }
\end{algorithm}
\vspace{-1cm}

\subsection{Equilibrating fluxes}
\label{sec:Eqlb_LocMinProb}
Flux equilibration requires the solution of a series of constrained minimisation problems \eqref{eq:eqlb_by_local_minimisation}.
This can be done directly \cite[Construction 3.4]{Ern_FluxEqlb_2015} or by splitting the process into an explicit part, followed by an unconstrained minimisation \cite[Appendix A]{Bertrand_HHO_2023}. 
Restricting this discussion on the second approach, the difference of equilibrated- and approximated flux is calculated in two steps:
\begin{equation}
    \fluxRHi - \hatfunc\fluxHi = \fluxRHziSEi + \fluxRHziSEii
    \label{eq:eqlb_se_split}
\end{equation}
While the determination of $\fluxRHzSEi$ is an interpolation-like task \cite[Algorithm 2]{Bertrand_HHO_2023}, $\fluxRHzSEii$ is determined on a patch-wise divergence free space \cite[Lemma 12]{Bertrand_HHO_2023}
\begin{equation*}
    \mathrm{V}^\Delta_m(\patch) := \left\{\vec{v} \in \feRT[m]:\; \div\vec{v}=0 \land \vec{v}\cdot\vec{n}=0 \;\text{on}\; \partial\patch\setminus\Gamma_\mathrm{D} \right\}\; \,
    \label{eq:patchwise_div_free_RT}
\end{equation*} 
based on the unconstrained minimisation problem
\begin{equation}
    \iprod[\patch]{\fluxRHzSEii}{\vec{v}_{z,\h}} = -\iprod[\patch]{\fluxRHzSEi}{\vec{v}_{z,\h}} \quad\text{for all}\quad \vec{v}_{z,\h} \in \mathrm{V}^\Delta_m(\patch) \; .
    \label{eq:eqlb_se_unconstrained_min}
\end{equation}
The therefore required hierarchic definition of the Raviart-Thomas space \cite[Appendix A.2]{Bertrand_HHO_2023} is implemented using Basixs \cite{Basix_2022} custom element.

From an algorithmic perspective, the semi-explicit equilibration is performed on each patch for multiple RHS simultaneously. 
This avoids re-assembly and re-factorisation (Cholesky decomposition) of the system matrix $\mat{A}$ of \eqref{eq:eqlb_se_unconstrained_min} for patches not attached to the boundary $\Gamma_\mathrm{N}$. 
The solution procedure is detailed in Algorithm \ref{alg:eqlb_se_patch_level}. 
\vspace{-0.5cm}
\begin{algorithm}
    \caption{Function: equilibrate\_flux\_semiexplt}
    \label{alg:eqlb_se_patch_level}

    \SetKwInOut{Input}{input }
    \SetKwFunction{cPatch}{Patch}\SetKwFunction{cpatch}{patch}

    \Input{$\left\{\fluxRHi\right\}$, $\left\{\Pi_{m-1}\,\fluxHi\right\}$, $\left\{\Pi_{m-1}\,\mathrm{f}_i\right\}$}

    \For{$i = 0;\ i < n_\mathrm{RHS};\ i++$}
    {
        Evaluate $\fluxRHziSEi$: \cite[Algorithm 2]{Bertrand_HHO_2023} using $\Pi_{m-1}\,\fluxHi$ and $\mathrm{f}_i$\;
        \BlankLine
        \uIf{$i = 0$}
        {
            Assemble $\mat{A}$ and $\vec{L}$ simultaneously and factorise $\mat{A}$\;
        }
        \Else
        {
            \lIf{$\partial\patch \cap \Gamma_\mathrm{N} = \emptyset $}
            {
                Re-assemble $\vec{L}$
            }
            \lElse
            {
                Assemble $\mat{A}$ and $\vec{L}$ simultaneously and factorise $\mat{A}$
            } 
        }
        \BlankLine
        Evaluate: $\mat{A}\cdot\fluxRHziSEii = \vec{L}$ and append solution: $\fluxRHi \mathrel{+}= \fluxRHziSEi + \fluxRHziSEii$\;
    }
\end{algorithm}
\vspace{-1cm}

\subsection{Imposing the weak symmetry}
\label{sec:Eqlb_WeakSym}
The weak symmetry condition is optionally enforced by an additional correction term resulting from a constrained minimisation problem. 
Performing flux equilibration on the rows of the stress tensor initially, yields a function fulfilling the required divergence- and boundary conditions.
Adding $\LinPKEqlbHzWS$ yields
\begin{equation}
    \LinPKEqlbHz - \hatfunc\LinPKEqlbH = \LinPKEqlbHzSEi + \LinPKEqlbHzSEii + \LinPKEqlbHzWS\; ,
\end{equation} 
fulfilling the weak symmetry condition.
Based on $\mat{J}(\bullet)=\begin{pmatrix} 0 & -1\\ 1 & 0 \end{pmatrix}\cdot (\bullet)$ and
\begin{equation*}
    \mathrm{V}_{1,0}(\patch) := \left\{\mathrm{v}_\h \in \HiSpace[\domain]:\, \mathrm{v}_\h\vert_\mathrm{T} \in \feP[k](\mathrm{T}) \land \iprod[\patch]{\mathrm{v}}{1} = 0 \;\text{if}\; \partial\patch \cap \Gamma_\mathrm{D} = \emptyset
    \right\}\; ,
\end{equation*}
the constrained solution $\left( \LinPKEqlbHzWS,\,\xi_{z,\h}\right) \in \mathrm{V}^\Delta_m(\patch) \times \mathrm{V}_{1,0}(\patch)$ satisfies
\begin{equation}
    \begin{split}
        \iprod[\patch]{\LinPKEqlbHzWS}{\bm{\tau}_{z,\h}} + \iprod[\patch]{\mat{J}\left(\xi_{z,\h}\right)}{\bm{\tau}_{z,\h}} &= 0\\
        \iprod[\patch]{\LinPKEqlbHzWS}{\mat{J}\left(\gamma_{z,\h}\right)} &= -\iprod[\patch]{\LinPKEqlbHzSEi + \LinPKEqlbHzSEii}{\mat{J}\left(\gamma_{z,\h}\right)}\; ,
    \end{split}
    \label{eq:eqlb_ws_constrmin}
\end{equation}
for all $\left(\tau_{z,\h},\,\gamma_{z,\h}\right) \in \mathrm{V}^\Delta_m(\patch)\times\mathrm{V}_{1,0}(\patch)$.
To guarantee solvability (see \cite{Bertrand_EqlbElast_2021}), patches must have at least two internal facets, and $k,m \geq 2$ has to hold.
In order to avoid the direct solution of \eqref{eq:eqs_ws_constrmin}, the equation system resulting from the saddle-point problem \eqref{eq:eqlb_ws_constrmin}, a Schur complement based solver is implemented ($\vec{u}_{R1}$ and $\vec{u}_{R2}$ are the DOFs for the rows of the stress tensor).
\begin{equation}
    \begin{bmatrix}
        \mat{A} & \mat{0} & \mat{B}_1 & \mat{0} \\
        \mat{0} & \mat{A} & \mat{B}_2 & \mat{0} \\
        \transpT{\mat{B}_1} & \transpT{\mat{B}_2} & \mat{0} & \vec{C}\\
        \vec{0} & \vec{0} & \transpT{\vec{C}}  & \mat{0}
    \end{bmatrix}\cdot
    \begin{bmatrix}
        \vec{u}_{R1} \\
        \vec{u}_{R2} \\
        \vec{c}\\
        \lambda
    \end{bmatrix}
    =
    \begin{bmatrix}
        \vec{0}\\
        \vec{0}\\
        \vec{L}_c\\
        0
    \end{bmatrix}
    \label{eq:eqs_ws_constrmin}
\end{equation}

\noindent Reusing the already factorised matrix $\mat{A}$ from flux equilibration the Schur complement $\mat{S} = \mat{B}_1\,\mat{A}^{-1}\,\transpT{\mat{B}_1} + \mat{B}_2\,\mat{A}^{-1}\,\transpT{\mat{B}_2}$ is calculated. 
Solving $\mat{S}\cdot\vec{c}=\vec{L}_c$ based on its LU decomposition followed by solving $\mat{A}\cdot\vec{u}_{Ri}=-\mat{B}_i\cdot\vec{c}$ gives a notable speed-up compared to a direct solution of \eqref{eq:eqs_ws_constrmin}. 
The general algorithmic structure is outlined in Algorithm \ref{alg:eqlb_ws_patch_level}.
\vspace{-0.5cm}
\begin{algorithm}
    \caption{Function: impose\_weak\_symmetry}
    \label{alg:eqlb_ws_patch_level}
    \SetKwInOut{Input}{input }
    \SetKwFunction{cPatch}{Patch}\SetKwFunction{cpatch}{patch}

    \Input{$\left\{\fluxRHi\right\}$}

    \BlankLine
    \lIf{$\left(\partial\patch \cap \Gamma_\mathrm{N} \neq \emptyset\right)$}
    {
        Assemble $\mat{A}$, $\mat{B}_i$, $\mat{C}$ and $\vec{L}_c$ simultaneously
    }
    \lElse
    {
        Assemble $\mat{B}_i$, $\mat{C}$ and $\vec{L}_c$ simultaneously
    } 

    \BlankLine
    \For{$i = 0;\ i < 2;\ i++$}
    {
        \lIf{$\left(\partial\patch \cap \Gamma_\mathrm{N} \neq \emptyset\right)$}
        {
            Apply BCs and re-factorise $\mat{A}$
        }
        $\mat{S} \mathrel{+}= \transpT{\mat{B}_i}\inv{\mat{A}}\mat{B}_i$\;
    }
    \BlankLine
    Solve: $\mat{S}\cdot\vec{c} = \vec{L}_c$\;
    \BlankLine
    \For{$i = 0;\ i < 2;\ i++$}
    {
        \lIf{$\left(\partial\patch \cap \Gamma_\mathrm{N} \neq \emptyset\right)$}
        {
            Apply BCs and re-factorise $\mat{A}$
        }
        Evaluate: $\mat{A}\cdot\vec{u}_{Ri} = -\mat{B}_i\cdot\vec{c}$\;
    }
    \BlankLine
    Append solution: $\LinPKEqlbHz \mathrel{+}= \LinPKEqlbHzWS$\;
\end{algorithm}
\vspace{-1cm}

\section{Results}
\label{sec:Results}
To illustrate the capabilities of dolfinx\_eqlb, adaptive solution procedures are presented for two characteristic problems.
Primal problems are solved using Lagrangian finite elements of degree $k$.
The performance of an error estimate $\eta$ is characterised based on the efficiency index $ \mathrm{i_{eff}} = \eta/\mathrm{err}$, where $\mathrm{err}$ denotes the true error .

\subsection{The Poisson equation}
As a first example, a rectangular domain with discontinuous coefficients $\kappa$ in each of the four quadrants is considered.
Two sets of parameters $\kappa_2=\kappa_4=1$ and $\kappa_1=\kappa_3$ with either $\kappa_1=5$ (20 refinement levels) or $\kappa_1=100$ (40 refinement levels) are tested.
Dirichlet BCs on $\Gamma_\mathrm{D}=\partial\Omega$ are prescribed based on the analytical solution $\GenSol_\mathrm{ext}$ from \cite{Riviere_ApostPoissonCoeffDG_2003}.
Meshes are refined based on a Dörfler marking strategy with $\theta=0.5$. 

\begin{figure}
    \centering
    \begin{subfigure}[b]{0.32\textwidth}
        \centering
        \begin{tabular}{@{}c|c c|c|c@{}}
            \toprule
            \multicolumn{1}{l}{$\kappa_1$} & $k$ & $m$ & $\mathrm{e.o.c}$ & $\mathrm{i_{eff}}$ \\ \midrule
            \multirow{4}{*}{$5$}     & 1 & 1 & $0.50$ & $1.47$ \\
                                     & 1 & 2 & $0.50$ & $1.06$ \\
                                     & 2 & 2 & $0.99$ & $1.40$ \\
                                     & 2 & 3 & $1.03$ & $1.05$ \\ \cmidrule(l){1-5} 
            \multirow{4}{*}{$100$}   & 1 & 1 & $0.50$ & $1.70$ \\
                                     & 1 & 2 & $0.54$ & $1.26$ \\
                                     & 2 & 2 & $1.36$ & $1.78$ \\
                                     & 2 & 3 & $2.14$ & $1.36$ \\ \cmidrule(l){1-5} 
        \end{tabular}
        \caption{}
        \label{fig:poisson-riviere_results}
    \end{subfigure}
    \hfill
    \begin{subfigure}[b]{0.32\textwidth}
        \centering
        \includegraphics[width=\textwidth]{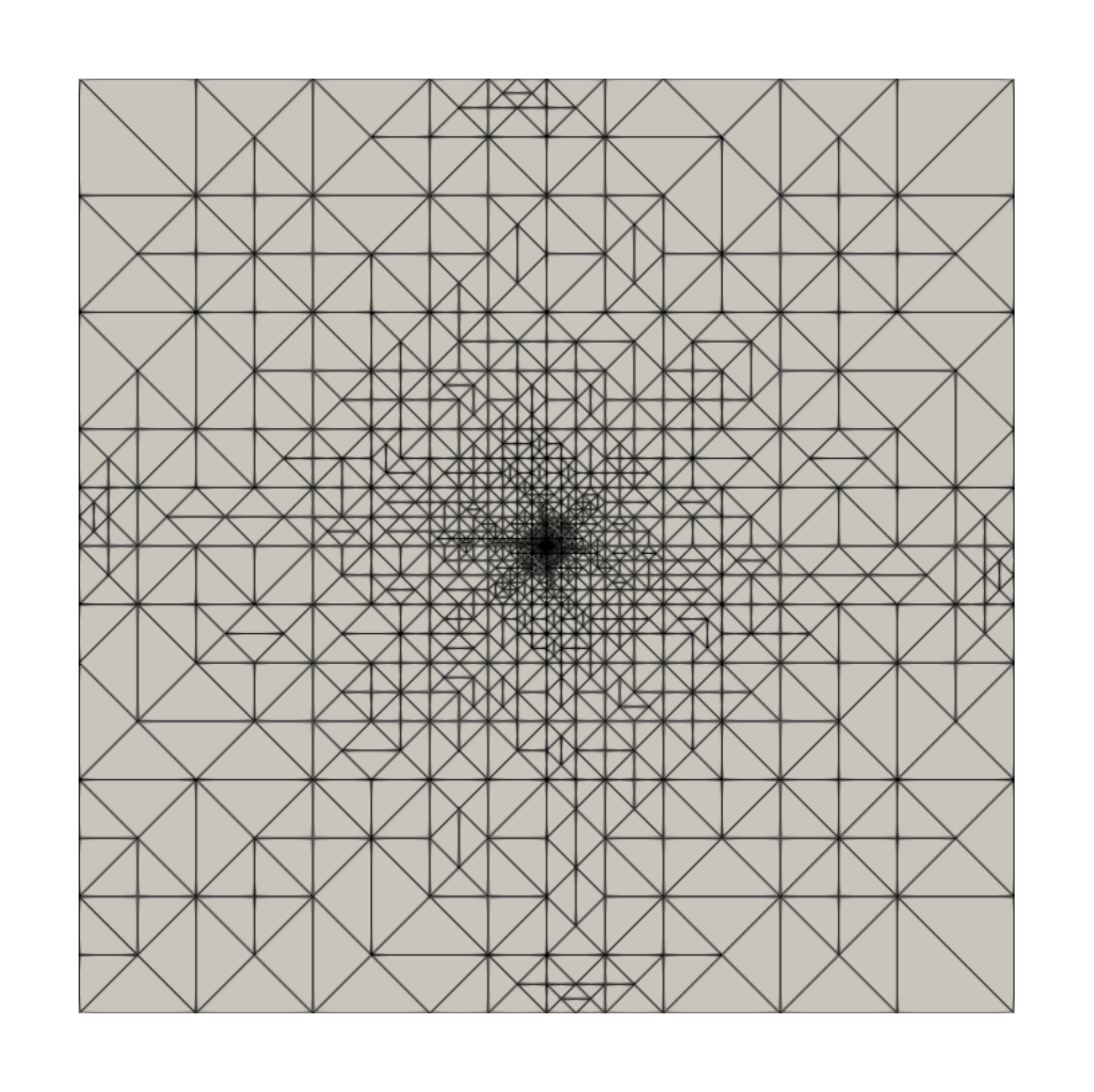}
        \caption{}
        \label{fig:poisson-riviere_refmesh-p1}
    \end{subfigure}
    \hfill
    \begin{subfigure}[b]{0.32\textwidth}
        \centering
        \includegraphics[width=\textwidth]{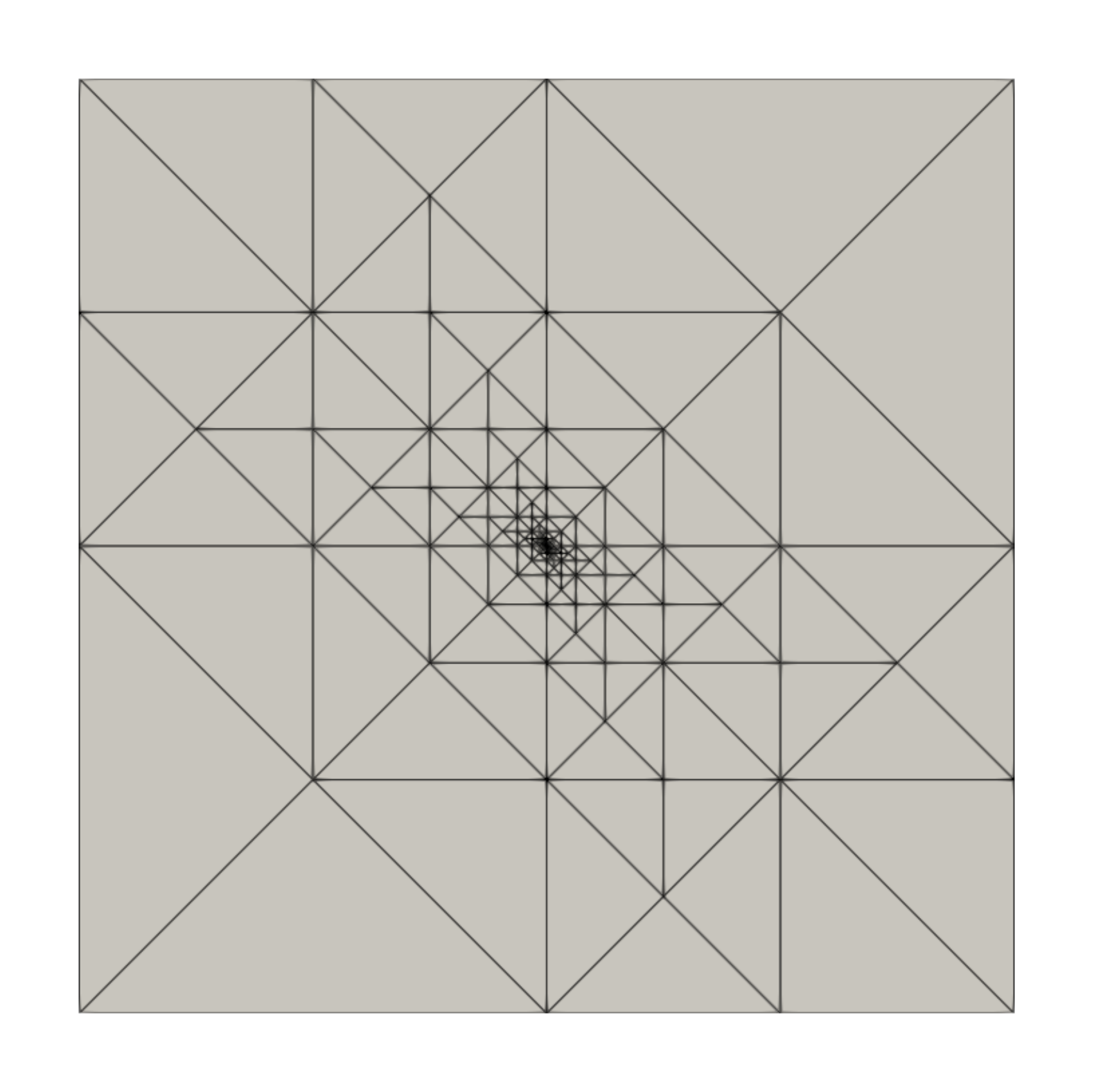}
        \caption{}
        \label{fig:poisson-riviere_refmesh-p2}
    \end{subfigure}
    \caption{Results of adaptive FEM calculations with different orders $k$ and $m$. E.o.c and $\mathrm{i_{eff}}$ after the final refinement step are reported in (a). The two final meshes for $\kappa_1=100$ are depicted in (b) for $k=m-1=1$ and (c) for $k=m-1=2$.}
    \label{fig:poisson-riviere}
\end{figure}

Convergence orders (e.o.c) and efficiency indices after the final refinement step are reported in Fig. \ref{fig:poisson-riviere_results}, while the final meshes for first- and second-order approximations of the case $\kappa_1=100$ are shown in Fig. \ref{fig:poisson-riviere_refmesh-p1} and \ref{fig:poisson-riviere_refmesh-p2}. 
The solutions are in good agreement with expectations from literature. 
Meshes are refined around the singularity in the center of the domain and the convergence rates are $\approx-0.5$ for $\GenSolH\in\feP[1]$ and $\approx-1$ for $\GenSolH\in\feP[2]$.
For $\kappa_1 = 100$  a convergence rate $>1$ indicates a pre-asymptotic state of convergence.
As further refinement would lead to cells with a characteristic length close to machine precision, the refinement has nevertheless been terminated.
The efficiency of the error estimate \eqref{eq:ee_poisson-disc-kappa} depends on the degree $m$ of the equilibrated flux.
Choosing $m=k+1$ yields efficiency indices close to one for $\kappa_1=5$ and between 1.2 and 1.4 for $\kappa_1=100$. 
For $m=k$ they are slightly worse.

\subsection{Linear elasticity}
\begin{figure}
    \centering
    \begin{subfigure}[b]{0.4\textwidth}
        \begin{tabular}{@{}l|c|c|ccc|l@{}}
            \toprule
            $k$ & $m$ & $n_\mathrm{DOF}$ & $\mathrm{err}$ & $\eta$ & $\eta _\mathrm{as}$ & $\mathrm{i_{eff}}$ \\ \midrule
            2 & 2 & 34070 & $0.0009$ & $0.009$ & $0.009$ & $10.7$ \\
            2 & 3 & 23202 & $0.0010$ & $0.008$ & $0.007$ & $7.9$ \\
            3 & 3 & 6656  & $0.0009$ & $0.020$ & $0.016$ & $17.0$ \\
            3 & 4 & 7100  & $0.0007$ & $0.009$ & $0.009$ & $13.0$ \\ \midrule
            2 & $2^*$ & 27788 & $0.0008$ & $0.001$ & $-$ & $1.5$ \\
            2 & $3^*$ & 26538 & $0.0008$ & $0.001$ & $-$ & $1.2$ \\
            3 & $3^*$ & 5738  & $0.0009$ & $0.001$ & $-$ & $1.5$ \\
            3 & $4^*$ & 5996  & $0.0008$ & $0.001$ & $-$ & $1.2$ \\ \bottomrule
        \end{tabular}
        \vspace{1.1cm}
        \caption{}
        \label{fig:cook_results-summary}
    \end{subfigure}
    \hfill
    \begin{subfigure}[b]{0.58\textwidth}
        \centering
        \includegraphics[scale=1.0]{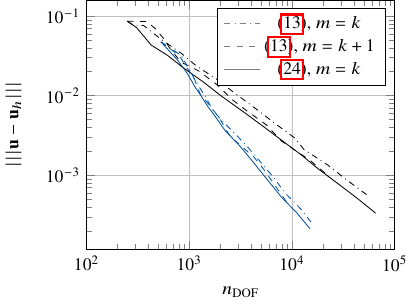}
        \caption{}
        \label{fig:cook_results-convhist}
    \end{subfigure}
    \caption{Effectivity of the different adaptive solution procedures for the Cooks membrane: (a) summaries the results for the first mesh with $\mathrm{err} = \vert\vert\vert \Displacement - \Displacement_\h \vert\vert\vert \leq 10^{-3}$, (b) details the convergence history (black: $k=2$, blue: $k=3$). Orders $m^*$ indicate the usage of \eqref{eq:elasticity_heuristic-ei}.}
    \label{fig:linelast-cook_results}
\end{figure}

\begin{minipage}[t]{0.65\textwidth}
Based on the Poisson equation, the influence of the equilibration order $m$ on the efficiency of the resulting estimate is shown.
When it comes to elasticity, the major difference is the weakly considered symmetry of stress tensor.
Influences of equilibration order on the error estimate are discussed followingly based on the Cooks membrane (see Fig. \ref{fig:cook_definition}), with $\tilde{\Lamei}=2.333$ and $t=0.03$. The mesh adaptivity is based on a Dörfler marking strategy with $\theta=0.6$.
Characteristics of the first mesh on which $\vert\vert\vert \Displacement - \Displacement_\h \vert\vert\vert \leq 10^{-3}$ holds, are summarised in Fig. \ref{fig:cook_results-summary}. As $\LinPKEqlbH$ exactly fulfils the divergence
\end{minipage}
\hfill
\begin{minipage}[t]{0.32\textwidth}
\centering\raisebox{\dimexpr-\height+\ht\strutbox}{
\includegraphics[scale=1.0]{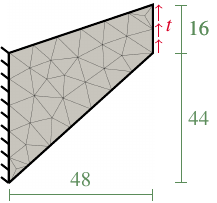}
}
\captionof{figure}{Cooks membrane.}
\label{fig:cook_definition}
\end{minipage}

\noindent condition from Definition \ref{def:equilibrated_stress}, $\eta$ reduced to the sum of the $\mathcal{A}$-norm of the stress difference $\LinPKEqlbH - \stressPKLin_\h$ and the norm of assymetric part of $\LinPKEqlbH$. 
Following \ref{fig:cook_results-summary}, the estimate is clearly dominated by the second part.
Using equilibrated stresses of order $m=k+1$ increases the efficiency, with a relative reduction (compared to the case with $m=k$) comparable to those from the Poisson example.
An increased accuracy of the error estimate affects the effectivity of the adaptive solution procedure -- measured by the number of degrees of
freedom, required for a certain error -- in a positive way.
This trend is clearly much more pronounced for $k=2$, whereby a similar accuracy is achieved with $47\%$ fewer degrees of freedom.

\newpage
\noindent A practical shortcut -- equilibration for $m=k+1$ is significantly more expensive than for $m=k$ -- seems to be the heuristic error indicator
\begin{eqnarray}
    \eta = \norm{\LinPKEqlbH - \stressPKLin_\h}\; ,
    \label{eq:elasticity_heuristic-ei}
\end{eqnarray}
where no weak symmetry is enforced on $\LinPKEqlbH$.
This yields efficiency indices close to one (see the bottom rows of Fig. \ref{fig:cook_results-summary}, and shows, comparing the convergence history in Fig. \ref{fig:cook_results-convhist}), slightly better results as the guaranteed estimate with $m=k+1$.

Up to this point only the accuracy of the error estimates has been compared. 
Within the following the focus will be on the total solution time $t_\mathrm{tot}$, which is made up from $t_\mathrm{prime}$ -- the time for assembly and solution, using PETSc and MUMPS, of the primal problem -- and $t_\mathrm{eqlb}$, the time for performing the equilibration.
Comparing in a first step the relative equilibration costs $t_\mathrm{eqlb}/t_\mathrm{tot}$ for the Cooks membrane with fixed meshes in Fig. \ref{fig:cook_performance-uniform}, $t_\mathrm{eqlb}$ is -- for a sufficient size of the primal problem -- considerably smaller than $t_\mathrm{prime}$.  
Increasing $k$ increases the effort required for equilibration, while using an equilibration with $m=k+1$ is significantly more expensive than the respective lowest order case $m=k$. Similar timings for the adaptive solution of the Cooks membrane (timings are accumulated until $\vert\vert\vert \Displacement - \Displacement_\h \vert\vert\vert \leq 10^{-3}$) are summarised in Fig. \ref{fig:cook_performance-adaptive}.
While for the cases $k=m=2$ the total solution time is dominated by the solution of the primal problem, this trend is reversed for more accurate (guarantied) estimate with $k=m-1=2$.
Even though having the fewest primal degrees of freedom, the entire solution time is the highest due to the high computational costs for the equilibration.
Using primal approximations based on $k=3$ reduces the overall computation time, but leading to a significant share of the equilibration in the total computation time. 
An effect, amplified by the small sizes of the primal problems.
As for $k=2$ the heuristic indicator \eqref{eq:elasticity_heuristic-ei} with $m=k$ performs best, while the higher order estimate with $m=k+1$ is outperformed.
Clearly these results will have to be reevaluated in a parallel context -- which is beyond the current scope -- are also for primal problems of larger size.
\begin{figure}[]
\centering
    \begin{subfigure}[b]{0.6\textwidth}
        \centering
        \begin{tabular}{@{}c|ccc||c|ccc@{}}
        \toprule
        \multicolumn{4}{c}{$k=2$}    & \multicolumn{4}{c}{$k=3$}    \\ \midrule
        $n_\mathrm{DOF}\,\setminus\, m$ & 2    & $2^*$    & 3    & $n_\mathrm{DOF} \,\setminus\, m$ & 3    & $3^*$    & 4    \\ \midrule
        $3.49 \cdot 10^3$  & 36.3 & 27.4 & 66.4 & $3.31 \cdot 10^3$ & 45.4 & 37.9 & 63.1 \\
        $1.36 \cdot 10^4$  & 22.0 & 15.1 & 49.4 & $1.30 \cdot 10^4$ & 33.6 & 27.0 & 51.7 \\
        $2.14 \cdot 10^5$  & 15.3 & 10.2 & 38.6 & $2.04 \cdot 10^5$ & 24.5 & 19.3 & 40.9 \\
        $8.54 \cdot 10^5$  & 12.6 & 8.34 & 33.7 & $8.13 \cdot 10^5$ & 20.7 & 16.2 & 35.6 \\ \bottomrule
        \end{tabular}
        \vspace{0.2cm}
        \caption{}
        \label{fig:cook_performance-uniform}
    \end{subfigure}
    \hfill
    \begin{subfigure}[b]{0.38\textwidth}
        \centering
        \begin{tabular}{@{}l|c|ccc@{}}
            \toprule
            $k$ & $m$ & $t_\mathrm{prime}$ [s] & $t_\mathrm{tot}$ [s] & ratio [\%] \\ \midrule
            2 & 2     & $0.47$ & $0.61$ & 23.2\\
            2 & $2^*$ & $0.38$ & $0.45$ & 16.4\\
            2 & 3     & $0.31$ & $0.65$ & 52.7 \\ \midrule
            3 & 3     & $0.09$ & $0.15$ & 41.9\\
            3 & $3^*$ & $0.07$ & $0.11$ & 35.6\\
            3 & 4     & $0.10$ & $0.26$ & 60.5\\ \bottomrule
        \end{tabular}
        \caption{}
        \label{fig:cook_performance-adaptive}
    \end{subfigure}
    \caption{Performance measurements based on the Cooks membrane. (a) $\mathrm{ratio} = t_\mathrm{eqlb} / t_\mathrm{tot}$ for different primal problems. (b) Accumulated timings using an adaptive algorithm until $\vert\vert\vert \Displacement - \Displacement_\h \vert\vert\vert \leq 10^{-3}$. Orders $m^*$ indicate the use of \eqref{eq:elasticity_heuristic-ei}.}
    \label{fig:linelast-cook_performance}
\end{figure}
\vspace{-0.8cm}

\section{Conclusions}
\label{sec:Conclusions}
Within this contribution dolfinx\_eqlb, a FEniCSx based library for the efficient equilibration of fluxes and stresses, has been introduced. 
Throughout characteristic examples for the Poisson problem and linear elasticity the efficiency of the resulting error estimates were shown.
An computationally efficient, but heuristic error indicator for elasticity is introduced, neglecting the asymmetry of the equilibrated stress.
Within our future work we intend to quantify computational efficiency of the presented implementation on problems of more realistic sizes. 
We further intend to generalise the implementation for 3D domains and transient problems like e.g. poroelasticity.

\subsection*{Supplementary material}
This work is based on dolfinx\_eqlb v1.2.0 (\url{https://github.com/brodbeck-m/dolfinx_eqlb/tree/v1.2.0}). 
The presented examples can either be accessed via GitHub (\url{https://github.com/ brodbeck-m/AFEM-by-Equilibration}) or, containing a Docker image, DaRUS (\url{https://doi.org/10.18419/darus-4500}).



\begin{thebibliography}{99.}%

\bibitem{FEniCSx_2023} Baratta, I.A. et al. (2023). DOLFINx: The next generation FEniCS problem solving environment. Zenodo. doi:10.5281/zenodo.10447666

\bibitem{Bertrand_Hypercircle_2020} Bertrand, F. and Boffi, D. (2020). The Prager–Synge theorem in reconstruction based a posteriori error estimation. Contemporary
Mathematics, 45–67. doi:10.1090/conm/754/15152

\bibitem{Bertrand_EqlbElast_2021} Bertrand, F., Kober, B., Moldenhauer, M. and Starke, G.: Weakly symmetric stress equilibration and a posteriori error estimation for linear elasticity. Comput. Math. Appl. (2021) doi: 10.1002/num.22741

\bibitem{Bertrand_HHO_2023} Bertrand, F., Carstensen, C., Gräßle, B. and Tran, N.T.: Stabilization-free HHO a posteriori error control. Numer. Math. (2023) doi: 10.1007/s00211-023-01366-8

\bibitem{Braess_EqlbFluxes_2008} Braess, D. and Schöberl, J.: Equilibrated Residual Error Estimator for Edge Elements. Math. Comput. \textbf{77}, 651--672 (2008)

\bibitem{Cai_SemiexplzEqlb_2012} Cai, Z. and Zhang, S.: Robust Equilibrated Residual Error Estimator for Diffusion Problems: Conforming Elements. SIAM J. Numer. Anal. (2012) doi: 10.1137/100803857

\bibitem{Ern_FluxEqlb_2015} Ern, A and Vohralík, M.: Polynomial-Degree-Robust A Posteriori Estimates in a Unified Setting for Conforming, Nonconforming, Discontinuous Galerkin, and Mixed Discretizations. SIAM J. Numer. Anal. (2015) doi: 10.1137/130950100

\bibitem{Prager_Equilibartion_1947} Prager, W. and Synge, J.L.: Approximations in elasticity based on the concept of function space. Q. J. Mech. Appl. Math. \textbf{5}, 241--269 (1947)

\bibitem{Riviere_ApostPoissonCoeffDG_2003} Rivière, B. and Wheeler, M.F.: A Posteriori error estimates for a discontinuous galerkin method applied to elliptic problems. Comput. Math. Appl. (2003) doi: 10.1016/S0898-1221(03)90086

\bibitem{Basix_2022} Scroggs, M.W., Dokken, J.S., Richardson, C.N. and Wells, G.N.: Construction of Arbitrary Order Finite Element Degree-of-Freedom Maps on Polygonal and Polyhedral Cell Meshes. ACM Trans. Math. Softw. (2022) doi: 10.1145/3524456

\end{thebibliography}
\end{document}